\documentclass[12pt]{amsart}
\usepackage{graphics}
\def\ints{\hbox{\sl Z\kern-.4em Z \kern-.3em}}
\def\inf{\mbox{\upshape inf}\,}
\def\sup{\mbox{\upshape sup}\,}
\def\cycl{\mathbf c}
\def\decy{\mathbf d}
\def\bQ{\mathbf Q}

\def\bD{\mathbf D}
\newtheorem{thm}{Theorem}
\newtheorem{lemm}[thm]{Lemma}
\newtheorem{cor}[thm]{Corollary}
\newtheorem{con}[thm]{Conjecture}

\begin{document}
\title[ braid conjugacy ]
{The infimum, supremum and geodesic length of a braid conjugacy class}
\author[Joan S.~Birman]{Joan S.~Birman$^1$\footnote{$^1$ Partially supported by
NSF Grants DMS-9705019 and DMS-9973232. }}
\address{Department of Mathematics, Barnard College of Columbia
University\\ 2990 Broadway\\ New York, NY 10027-4427}

\email{jb@math.columbia.edu}

\author[Ki Hyoung Ko]{Ki Hyoung Ko$^2$\footnote{$^2$ The author wishes to
acknowledge the financial support of the Korea Research Foundation made in the
program year of 1999}}
\author{Sang Jin Lee}
\address{Department of Mathematics\\
Korea Advanced Institute of Science and Technology\\
Taejon, 305--701, Korea}

\email{\{knot, sjlee\}\char`\@knot.kaist.ac.kr}

\begin{abstract}
\noindent Algorithmic solutions to the conjugacy problem in the braid groups
$B_n, n=2,3,4,\dots$ were given in \cite{EM1994} and in \cite{BKL1998}.  This note
concerns the computation of two integer class invariants, known as `inf' and
`sup'. A key issue in both algorithms is the number
$m$ of times one must `cycle' (resp. `decycle') in order to either increase inf
(resp. decrease sup) or to be sure that it is already maximal (resp. minimal) for
the class. Our main result is to prove that $m$ is bounded above by
$((n^2-n)/2)-1$ in the situation of \cite{EM1994} and by
$n-2$ in the situation of \cite{BKL1998}.  It follows immediately that
the computation of inf and sup is polynomial in both word length and braid
index, in both algorithms.   The integers inf and sup determine (but are not
determined by) the shortest geodesic length for elements in a conjugacy class,
as defined in \cite{Charney}, and so we also obtain a polynomial-time algorithm
for computing this length.
\end{abstract}

\maketitle
\centerline{March 21, 2000}

\section{{\bf Introduction}}
The conjugacy problem in the $n$-string braid group $B_n$  is
the following decision problem:

\begin{quote}
Given two braids $\alpha,\alpha^\prime \in B_n$,
determine, in a finite number of steps, whether
$\alpha=\gamma\alpha^\prime\gamma^{-1}$ for some $\gamma\in B_n$.
\end{quote}
In the late sixties Garside  ~\cite{Garside}
solved the (word and) conjugacy problems in $B_n$.  His solution
to both problems was exponential in both word length and braid index.
Subsequently, the efficiency of his algorithm was improved by
Thurston~\cite{Thurston1992} and Elrifai-Morton~\cite{EM1994} to give a solution
to the word problem which is polynomial in both word length and braid index.

All three papers
\cite{Garside}, ~\cite{Thurston1992} and ~\cite{EM1994} work with the
following well-known presentation of $B_n$, which we will call the \emph{old
presentation}:

\begin{quote}
$\begin{array}{ll}
\mbox{generators:}& \sigma_1,\ldots,\sigma_{n-1} \\
\mbox{relations:} & \sigma_i\sigma_j=\sigma_j\sigma_i,\quad |i-j| > 1\\
& \sigma_i\sigma_j\sigma_i = \sigma_j\sigma_i\sigma_j,\quad |i-j|=1
\end{array}$
\end{quote}
\noindent There is also a parallel and slightly more efficient solution to the
word and conjugacy problems in \cite{BKL1998},  due to the authors of this paper.
It uses  a different presentation which we call the  \emph{new presentation}:

\begin{quote}
$\begin{array}{ll}
\mbox{generators:}& a_{ts},\quad n\ge t>s\ge 1 \\
\mbox{relations:} & a_{ts}a_{rq}=a_{rq}a_{ts},\quad (t-r)(t-q)(s-r)(s-q)>0\\
 & a_{ts}a_{sr}=a_{tr}a_{ts}=a_{sr}a_{tr},\quad n\ge t>s>r\ge 1\end{array}$
\end{quote}
\noindent The terms {\em old} and {\em new} are due to Krammer, who used the new
presentation in \cite{Krammer}. Both the old and new solutions to the  word problem
are polynomial in word length and braid index, but the  best estimates
obtained for the complexity of the solution to the  conjugacy problem (see
\cite{BKL1998}) were rough exponential bounds.  It was clear that better  answers
could not be obtained without more detailed information about the  combinatorics,
using either the old or new presentation.

Let $|W|$ denote the letter length of $W$, as a word in the
given set of generators of $B_n$.  The main result in this note is an
algorithm which is polynomial in both $|W|$ and $n$ for computing two key integer
invariants of the conjugacy class $[W]$ of $W$. The invariants in question are
known as the {\it infimum} and {\it supremum} (or more informally {\it inf}
and {\it sup}), using either presentation. See $\S$2 below for precise
definitions. We will also be able to compute the {\em geodesic length} (defined
in $\S$2,4 below) for the conjugacy class in polynomial time.

The reason we are able to do this requires some explanation. The method for
finding {\it inf} (resp. {\it sup}) in both \cite{EM1994} and \cite{BKL1998}
rests on a procedure which is known as {\em cycling}  (resp. {\em decycling}).
While cycling and decycling are clearly finite  processes, it had not been
known how many times one must iterate them to  either increase inf($W'$) (resp.
decrease sup($W'$)) for a word $W'\in [W]$ or to guarantee that a
maximum (resp. minimum) value, denoted by inf($[W]$) (resp.
sup($[W]$)), 
for the conjugacy class has already been
achieved. For the old presentation it had been claimed in
~\cite{Thurston1992}
that the bound is 1, however an example was given in ~\cite{EM1994}
for which 2 cyclings were needed to increase the infimum.
Up to now, there were no published results which gave bounds, except for a
very crude estimate in \cite{BKL1998}.  Our main result in this note is
to  find upper and lower bounds for the number of times one must
cycle (resp. decycle), using either presentation,
 in order to replace  a given word $W$ with
$W^\prime\in[W]$, where $\inf(W^\prime)>\inf(W)$
(resp. $\sup(W^\prime)<\sup(W))$, or be sure that $W$
realizes inf($[W]$) (resp. sup($[W]$)).  For the new presentation we will
prove that our upper bound is the best possible one.

Here is an outline of this paper. In $\S$\ref{sec:statement of results} we
review the background and state our results in a precise way. See
Theorem \ref{thm:cyclbd}, Corollary \ref{cor:complexity estimates} and
Corollary \ref{cor:geodesic length}.  In
\S\ref{sec:proof}, we prove these three results.  In \S\ref{sec:ex},
we give examples which prove that the bound in Corollary
\ref{cor:complexity estimates} is sharp for the new presentation, with somewhat
weaker results for the old. In
$\S$\ref{sec:complexity issues} we discuss the open problem of whether the
solutions which we know to the conjugacy problem are polynomial in word length and
braid index, and state several conjectures relating to that matter and also to the
`shortest word problem' in $B_n$, defined in that section.

\section{{\bf Statement of Results}}
\label{sec:statement of results}

In this section we state our results precisely.  To do so we need to review
what has already been done. Since almost all the machinery
is identical in the two  theories, it will be convenient to introduce
unified notation, so that we may review both theories at the same time. The
symbol $W$ will be used to indicate a word in the generators of $B_n$, using
either presentation. The element and conjugacy class which $W$ represents will
be denoted $\{W\}$ and $[W]$.  The letter length of $W$ is  $|W|$.

\smallskip

\begin{enumerate}

\item [2.1] Note that the relations in the old and new presentations are
equivalences between positive words with same word-length. So the word-length
is easy to compute for positive words.
Let $B_n^+$ be the semigroup defined by the same generators and relations
in the given presentation.
The natural map $B_n^+\to B_n$ is
injective.~\cite{Garside, BKL1998}.

\smallskip

\item [2.2] There is a \emph{fundamental braid} $\bD$.
In the old presentation, $\bD$ has length $((n^2-n)/2) - 1$ and is the
half-twist
$$\Delta=(\sigma_1\cdots\sigma_{n-1})(\sigma_1\cdots\sigma_{n-2})
  \cdots(\sigma_1\sigma_2)\sigma_1.$$
In the new presentation it has length $n-1$ and it is the $(1/n)$-twist
$$\delta=a_{n(n-1)}a_{(n-1)(n-2)}\cdots a_{32}a_{21}.$$
The fundamental braid admits many many braid transformations, in both the old
and the new presentations, and so can be written in many ways as a positive
word in the braid generators. As a result of this flexibility, it has two important
properties:

\smallskip

\begin{enumerate}
\item [(i)] For any generator $a$,
there exist $A,B\in B_n^+$ such that: \\
$\bD=aA=Ba$;
\item[(ii)] For each generator $a$ we have $a\bD = \bD\tau(a)$
and also $\bD a = \tau^{-1}(a)\bD$,
where $\tau$ is the automorphism of $B_n$ which is defined by
$\tau(\sigma_i)=\sigma_{n-i}$  for the old presentation
and $\tau(a_{ts})=a_{(t+1)(s+1)}$ for the new presentation.
\item [(iii)] $\tau(\{\bD\}) = \{\bD\}$.
\end{enumerate}

\smallskip

\item [2.3] There are partial orderings `$\ge$' and `$\le$' in $B_n$.
For two words $V$ and $W$ in $B_n$ we say that
$V\ge W$ (resp.\ $W\le V$) if $V=PW$ (resp.\ $V=WP$) for some $P\in
B_n^+$. Note that $W$ is a positive word if and only if
$W\ge e$. We denote $V<W$ (resp.\ $V>W$) if $V\le W$ (resp.\ $V\ge W$) and $V\ne
W$.  In general $V\ge W$ is not equivalent to $W\le V$, although if either
$W$ or $V$ is a power of $\bD$ the two ordering conditions are equivalent
because powers of $\bD$ commute with elements of $B_n$ up to powers of the
index-shift automorphism $\tau$. Note that $\tau$ preserves the partial
ordering.

\smallskip

\item [2.4] The symbol $\bQ$ denotes the set of all initial subwords of $\bD$,
and $\bQ^\star = \bQ\setminus \{e,\bD\}.$  The cardinality
$|\bQ_{old}|$ is $n!$, whereas the cardinality $|\bQ_{new}|$ is the $n^{th}$
Catalan number. Note that
$|\delta| < |\Delta |$, also $|\bQ_{new}| < |\bQ_{old}|$. These are
the main reasons why it is sometimes easier to work with the new presentation
than the old.

\smallskip

\item [2.5] The {\it geodesic length} $l_Q(\{W\})$
was introduced and investigated by Ruth Charney in \cite{Charney}. It is the
smallest integer $k$ such that there is a word $q_1q_2\cdots q_k$ representing
$\{W\}$, with each $q_i\in\bQ\cup\bQ^{-1}$.  Define the geodesic length of the
conjugacy class $l_Q([W])$ to be the shortest such representation for words in the
conjugacy class $[W]$.

\smallskip

\item [2.6]
For each positive word $P$, there is a decomposition, called the
\emph{left-greedy decomposition}, $P=A_0P_0$ for $A_0\in\bQ$
and $P_0\ge e$, where $A_0$ has maximal length among all such decompositions,
i.e.\ if $P=A'_0P'_0$, where $A'_0\in \bQ$ and $P'_0\in B_n^+$, then
$A_0'\le A_0$.  The term `greedy' suggests that $A_0$ has absorbed as many
letters from $P_0$ as it can without leaving $\bQ$. The canonical factor
$A_0$ is called the \emph{maximal head} of $P$.  If $P = A_0P_0 =
A_0'P_0'$ in left greedy form, then $\{A_0\}=\{A_0'\}$ and
$\{P_0\}=\{P_0'\}$.  (Remark: The term {\em
left-canonical decomposition} was used in
\cite{BKL1998} and \cite{EM1994}, however in
recent years {\it left-greedy decomposition} has become the term of choice for
the same concept in the literature, hence we now change our notation.)

\smallskip

\item [2.7]
Any word $W$ in the generators admits a unique {\it normal form} which solves the
word problem in $B_n$. The normal form is:
$$ W=\bD^u A_1A_2\cdots A_k,\qquad
   u\in\ints,\; A_i\in\bQ^\star,$$
where for each $1\le i\le k-1$, the product $A_iA_{i+1}$ is
a left-greedy decomposition.
The integer $u$ (resp.\ $u+k$) is called the \emph{infimum} of $W$
(resp.{\em supremum} of $W$)
and denoted by $\inf(W)$ (resp.\ $\sup(W)$).

\smallskip

\item [2.8]
To solve the conjugacy problem, we need to study the maximum and minimum values
of inf and sup for the conjugacy class rather than for the word class. We consider
the following two operations
$\cycl$ and $\decy$, called \emph{cycling} and \emph{decycling},
respectively. For a given braid in normal form
$W=\bD^u A_1A_2\cdots A_k$, we define:
\begin{eqnarray*}
\cycl(W) &=& \bD^u A_2A_3\cdots A_k\tau^{-u}(A_1)\\
\decy(W) &=& \bD^u \tau^u(A_k) A_1 \cdots A_{k-1}.
\end{eqnarray*}
In general the braids on the right hand side will {\em not} be in
normal form, and must be rearranged into normal form
before the operation can be repeated.

\smallskip

\item [2.9] {\bf Theorem}
(see ~\cite{BKL1998,EM1994}):
\label{thm:cycl}

\smallskip

\begin{enumerate}
\item[(1)] If $W$ is conjugate to $V$ and if $\inf(V)>\inf(W)$, then
repeated cycling
will produce $\cycl^\ell(W)$ with $\inf(\cycl^\ell(W))>\inf W$.
\item[(2)] If $W$ is conjugate to $V$ with $\sup(V)<\sup(W)$, then
repeated decycling
will produce $\decy^\ell(W)$ with $\sup(\decy^\ell(W))<\sup(W)$.

\item[(3)] The maximum value of $\inf$ and the minimum value of $\sup$
can be achieved simultaneously.
\end{enumerate}

\smallskip

\item [2.10] The \emph{super summit set} SSS($[W]$) (\cite{BKL1998,EM1994}) is
the set of all conjugates  of $W$ which have the maximal infimum and the minimal
supremum in the conjugacy class $[W]$. It is a proper subset of the {\em summit
set} SS$([W])$ which was introduced in by Garside in \cite{Garside}.

\smallskip

\item [2.11] {\bf Theorem} (see \cite{BKL1998, EM1994,
Thurston1992}):
\label{thm:convex}
Let $W, W'$ be any two words in SSS($[W]$ = SSS$([W'])$.
Then there is a sequence
\[ W=W_0 \to W_1\to \cdots \to W_k=W' \]
such that each intermediate braid $W_i\in {\rm SSS}([W])$ and each $W_{i+1}$ is a
conjugate of $W_i$ by a single member of $\bQ$.

\smallskip

\item [2.12] By the theorems in $\S$2.9 and $\S$2.11
one can compute SSS($[W]$) as follows:
\begin{itemize}
\item Obtain an element $W'$ in the super summit set
by iterating cyclings and decyclings, starting with any given word $W$.
\item Compute the whole super summit set from $W'$ as follows:
Compute $AW'A^{-1}$ for all $A\in \bQ$ and
collect the braids in the super summit set.
Repeat the same process with each newly obtained element,
until no new elements are obtained.
\end{itemize}
Therefore there is a finite time algorithm to generate SSS$(W)$.
This algorithm solves the conjugacy problem in $B_n$. The
integers $\inf([W])$ and $\sup([W])$ are the same for all members of SSS$(W)$
and so are partial invariants of the conjugacy class $[W]$.
\end{enumerate}

\smallskip

In this article, we obtain an upper bound for the necessary number of cyclings
and decyclings in the theorem in $\S$2.9 above.
for both the old presentation and the new presentation.
We denote the word length of $W$ by $|W|$. Our main result is:
\begin{thm}
\label{thm:cyclbd}
Let\/ $W\in B_n$.  If\/
$\inf(W)$ is not maximal for
$[W]$, then\/ \inf$(\cycl^{|\bD|-1}(W)) > \inf(W)$.  If\/ $\sup(W)$ is not
minimal for $[W]$, then \sup$(\decy^{|\bD|-1}(W)) < \sup(W)$.
\end{thm}

\smallskip

As immediate applications, we have:

\begin{cor}
\label{cor:complexity estimates} Given any braid word $W\in B_n$,
there is an algorithm which is polynomial in both word length and
braid index for the computation of $\inf[W]$ and $\sup[W]$. Using
the new presentation the complexity of the algorithm is
$O(|W|^2n^2)$. %%2
\end{cor}

\begin{cor}
\label{cor:geodesic length}
There is an algorithm which is polynomial in both word length and braid index
for the computation of the geodesic length $l_Q([W])$ of the conjugacy class of
$W$, using either presentation. Using the new presentation the complexity is
$O(|W|^2n^2)$. %%3
\end{cor}

\smallskip

\section{{\bf Proof of Theorem \ref{thm:cyclbd} and Corollaries
\ref{cor:complexity estimates} and \ref{cor:geodesic length}.}}
\label{sec:proof}

\noindent{\bf Proof of Theorem \ref{thm:cyclbd}:}
We focus on  cycling because the proof  and the
difficulties are essentially identical  for decycling.

Here is the plan of the proof. We begin with a word $W=\bD^u P$ which is in
normal form, so that $u=\inf(W)$ and $P>e$. By hypothesis $\inf([W]) > u$, so
there exists an integer
$m$ such that $u = \inf(W) = \inf(\cycl(W))=\cdots =  \inf(\cycl^m(W))$, but
$\inf(\cycl^{m+1}(W)) > u$.
Each instance of cycling can be
realized by conjugation of $W$ by an element in
$\bQ^\star$, so we know there are
$A'_1,A'_2\dots,A'_m\in\bQ^\star\cup(\bQ^\star)^{-1}$ %%4
such that after conjugating
$W$, successively, by $A'_1,A'_2\dots,A'_m$ we  obtain $W'= R'\bD^uP(R')^{-1}$
with $\inf(W')=u+1$. (See Lemma \ref{lemm:existence of R} below.) Write $R'$ in
normal form. (See Lemmas \ref{lemm:inf(R)} and \ref{lemm:poswd}.) Our plan is to
show that the sequence of lengths of the canonical factors
$H'_m,\dots,H'_0$ for $R'$ satisfies
$|H'_m| < |H'_{m-1}|<\dots <|H'_0|$. %%5
Since each $H'_i\in\bQ^\star$, we
have $e<|H'_i|<|\bD|$. This places a limit on the length of the chain, i.e.
$m+1\leq |\bD|-1$ or %%6
$m\leq |\bD|-2$, as claimed.

We used the symbols $R',A_i',H_j'$ in the description above, but in the
actual proof we will use symbols $R,A_i,H_j$ which differ a little bit
from $R',A_i',H_j'$ because we wish to focus on the changes in the positive
part $P$ of $W$, rather than on changes in $\bD^uP$:

\smallskip

\begin{lemm}
\label{lemm:existence of R}
Choose any $W\in B_n$. Let $W=\bD^uP$, where $u=\inf(W)$.  Then
$\inf([W])> \inf(W)$
 if and only if there exists a
positive word $R$ such that $R P \tau^{-u}(R^{-1})\ge\bD$.
\end{lemm}

\begin{proof}
By hypothesis $\inf([W]) > \inf(W)$, so there exists
$X\in B_n$ with $\inf(XWX^{-1}) > \inf(W).$  Let
 $X = \bD^v Y, Y\ge e$, where $ v = \inf(Y)$. Then:
$$ (\bD^vY)(\bD^uP)(Y^{-1}\bD^{-v})\ge\bD^{u+1},$$
which implies (via part (ii) of (3) above) that:
$$ (\tau^u(Y))(P)(Y^{-1}) \geq \bD.$$
Set $R = \tau^u(Y)$, so that $Y = \tau^{-u}(R)$. Then $R\ge e$ and
$$ (R)(P)(\tau^{-u}(R^{-1})) \ge \bD,$$ as claimed.
\end{proof}

\smallskip

\noindent We will need to understand the structure of the positive word
$R$ in Lemma
\ref{lemm:existence of R}, and to learn how the normal form of
$R$ is related to that of $W$ and its images under repeated cycling.
Once we understand all these issues, we will be able to extract information
from $R$ about repeated cycling.  We begin our work with several
preparatory lemmas (i.e. Lemmas \ref{thm:head}, \ref{lemm:|W|+1} and
\ref{thm:inv}.):

\smallskip

\begin{lemm}
\label{thm:head}
Suppose that $P\ge e$ and that $RP\ge \bD$ for some $R\ge e$.
Let $P=A_0P_0$ be in left-greedy form.
Then $RA_0\ge \bD$.
\end{lemm}

\begin{proof}
See  Proposition 3.9 (IV) of~\cite{BKL1998}
for the new presentation and Proposition 2.10
of~\cite{EM1994} for the old presentation.
\end{proof}

\smallskip

\begin{lemm}
\label{lemm:|W|+1}
If $W\in B_n$ and $A\in \bQ$, then $\inf(W)\le \inf(WA) \le \inf(W)+1.$
\end{lemm}

\begin{proof}
Since $W\le WA\le W\bD$ and $\inf(W\bD)= \inf(W)+1$ the
assertion follows.
\end{proof}

For each $A\in\bQ$, let $\bar A$ denote the unique member of $\bQ$
which satisfies $\bar A A=\bD$.

\begin{lemm}
\label{thm:inv}
Let $Z = B_\ell B_{\ell-1}\cdots B_1$ be the normal form for $Z\ge e$.
Then $Z^{-1}\bD^\ell\ge e$ and the normal form
for $Z^{-1}\bD^\ell$ is
$$\tau(\bar B_1)\tau^2(\bar B_2)\cdots \tau^\ell(\bar B_\ell).$$
\end{lemm}

\begin{proof}
Observe that $\bD = \bar B_i B_i$ implies that
$\bD = \tau^i(\bar B_i) \tau^i(B_i)$ for every $i=1,\dots,\ell$.  Therefore:
\begin{eqnarray*}
Z^{-1}\bD^\ell
&=& B_1^{-1}B_2^{-1}\cdots B_\ell^{-1}\bD^\ell\\
&=& (\bD\tau(B_1^{-1}))(\bD\tau^2(B_2^{-1}))\cdots(\bD\tau^\ell
(B_\ell^{-1}))\\
&=& (\tau(\bar B_1)\tau(B_1)\tau(B_1^{-1}))(\tau^2(\bar B_2)
\tau^2(B_2)\tau^2(B_2^{-1}))\cdots (\tau^\ell(\bar B_\ell)\tau^\ell(B_\ell)
\tau^\ell(B_\ell^{-1}))\\
&=&  \tau(\bar B_1)\tau^2(\bar B_2)\cdots\tau^\ell(\bar B_\ell)
\end{eqnarray*}
because
$$\tau(B_i)\tau(B_i^{-1}) = \tau(B_iB_i^{-1}) = \tau(e) = e.$$
\end{proof}

\smallskip

We continue the proof of Theorem \ref{thm:cyclbd}. By $\S$2.8
and Lemma \ref{lemm:|W|+1}
there exists a nonnegative integer $m$
such that $\inf(\cycl(W)) = \inf(\cycl^2(W)) = \cdots = \inf(\cycl^m(W))$ but
$\inf(\cycl^{m+1}(W)) = \inf(W) + 1$.  To prove
Theorem~\ref{thm:cyclbd}, we must show that
$m+1$ is bounded above by $|\bD|-1$. Let $W = \bD^u P$, where $P\ge e$ and $u =
\inf(W)$.  By Lemma \ref{lemm:existence of R}
we know that $\inf([W])> \inf(W)$
 if and only if there exists a
positive word $R$ such that $R P \tau^{-u}(R^{-1})\ge\bD$.
Assume that among all such words we have chosen $R$ so that
$|R|$ is minimal.
We wish to describe this shortest word $R$ as a specific product
(in general not left-greedy) of canonical factors.
Our first observation is:

\begin{lemm}
\label{lemm:inf(R)}
$\inf(R)=0$.
\end{lemm}

\begin{proof} If not, then
$R=\bD R'$ for some $R'\ge e$ and
 $$R' P\tau^{-u}(R')^{-1} = \bD^{-1} RP\tau^{-u}(R)^{-1}\bD =
\tau(RP\tau^{-u}(R)^{-1})\ge \bD,$$ which contradicts the
minimality of $|R|$.
\end{proof}

\smallskip

\begin{lemm}
\label{lemm:poswd}
Let $\cycl^i(W)=
\bD^u A_iP_i$, where $A_i$ is the maximal head of $\cycl^i(W)$.
Then the positive word $R$ whose existence is guaranteed by
Lemma \ref{lemm:existence of R} is related to the $A_i$'s as follows:
$$R=\tau^{-m}(\bar A_m)\cdots \tau^{-1}(\bar A_1) \bar A_0. $$
\end{lemm}

\begin{proof}
Our starting point is:
$$R P \tau^{-u}(R)^{-1}\ge\bD, $$
which implies that $RP\ge
\bD$. Since
$P=A_0P_0$ is left-greedy, Lemma \ref{thm:head} then implies that
$RA_0\ge\bD$  and so $R= R_1 \bar A_0$ for some positive
word $R_1$. Now:
\begin{eqnarray*}
RA_0P_0\tau^{-u}(R)^{-1}
&=& R_1 \bar A_0 A_0 P_0 \;\;
    \tau^{-u}(\bar A_0^{-1}) \tau^{-u}(R_1^{-1})\\
&=& R_1 \bD  P_0 \tau^{-u}(A_0) \bD^{-1}
\tau^{-u}(R_1^{-1})\\
&=& R_1 \tau^{-1}\left( A_1P_1 \right) \tau^{-u}(R_1^{-1}).
\end{eqnarray*}

Since $RA_0P_0\tau^{-u}(R^{-1})\ge\bD$, we conclude that:
%% $$\tau(R_1) A_1P_1\tau^{1-u}(R_1^{-1}) \ge \bD.$$
%% and so also that:
$$R_1\tau^{-1}(A_1P_1)\tau^{-u}(R_1^{-1}) \ge \bD.$$  %%7
Iterating the construction, we obtain $R_1=R_2\tau^{-1}(\bar A_1)$ for
some positive word $R_2$, also $R_2 = R_3\tau^{-2}(\bar A_2),\dots,
R_m = R_{m+1}\tau^{-m}(\bar A_m)$.  Putting all of these together we learn that:
$$R = R_{m+1}\tau^{-m}(\bar A_m)\cdots \tau^{-1}(\bar A_1) \bar A_0 $$
for some positive word $R_{m+1}$.
Let $S = \tau^{-m}(\bar A_m)\cdots \tau^{-1}(\bar A_1) \bar A_0,$ so that $R =
R_{m+1}S.$  A straightforward calculation shows that %%8
$$\tau^{-(m+1)}(SP\tau^{-u}(S^{-1}))= A_{m+1}P_{m+1}.$$
Since
$\inf(\cycl^{m+1}(W)) = \inf (W)+1$, we have
$$1=\inf(\tau^{-(m+1)}(SP\tau^{-u}(S^{-1})))
=\inf(SP\tau^{-u}(S^{-1})).$$
By the minimality of $|R|$, we must have $R=S$.  Lemma
\ref{lemm:poswd}  is proved.
\end{proof}

\smallskip

The expression given for $R$ in the statement of Lemma \ref{lemm:poswd}
is in general not in normal form. We now study the maximal head $H_0$
of $R = \tau^{-m}(\bar A_m)\cdots \tau^{-1}(\bar A_1) \bar A_0$, and related
canonical factors $H_1,\dots,H_m$.   To define them, let $H_k$ be the maximal head of
$\tau^{-m}(\bar A_m) \cdots\tau^{-k}(\bar A_k)\subset R.$

\begin{lemm}\label{thm:headne}
$ e<H_m<H_{m-1}<\cdots <H_1< H_0<\bD. $
\end{lemm}

\begin{proof}
Our first observation is that $\inf(R)=0$ (see Lemma
\ref{lemm:inf(R)}). Since $H_0$ is the maximal head of $R$, it
follows that $$H_0<\bD.$$ Our second observation is that by
hypothesis $\inf(\cycl^m(W))=u$ and $\cycl^m(W) = \bD^u A_mP_m$ is
left-greedy, so that $A_m<\bD$, which implies that $e<\bar A_m$
and so: $$e<H_m.$$ Our third observation is that by the definition
of $H_k$ we must have: $$ H_m \le H_{m-1}\le \cdots \le H_1\le
H_0. $$ Therefore the only thing that we need to prove is that
$H_{k+1}\ne H_k$ for $k=0,1,\ldots, m-1$.\\

We first prove the assertion for $k=0$.
Assume that $H_1=H_0$. We will show that this leads to  a contradiction
to our choice of $R$.\\

We are given that:
$$R= \tau^{-m}(\bar A_m)\cdots \tau^{-1}(\bar A_1) \bar A_0 =
B_\ell B_{\ell-1}\cdots B_1,$$
where the decomposition on the left comes from Lemma \ref{lemm:poswd}
and the one on the right is the normal form for $R$.
By Lemma \ref{thm:inv} the normal form for $R^{-1}\bD^\ell$ is
$\tau(\bar B_1) \tau^2(\bar B_2) \cdots \tau^\ell(\bar B_\ell).$

By hypothesis  $H_1 = H_0 = B_\ell$, so
$$\tau^{-m}(\bar A_m) \cdots \tau^{-1}(\bar A_1)= B_\ell R_1$$
for some $R_1\ge e$.
Since $B_\ell R_1 \bar A_0 = R = B_\ell B_{\ell-1}\cdots B_1$,
it follows that $B_{\ell-1}\cdots B_1 = R_1\bar A_0$ and so
$$B_{\ell-1}\cdots B_1 P\ge \bD.$$
Let $a_i$ be the infimum of
$B_{\ell-1}\cdots B_1 P\tau^{-u}\left(
  \tau(\bar B_1) \tau^2(\bar B_2) \cdots \tau^i(\bar B_i) \right)$.
Then
\begin{enumerate}
\item $a_0\ge 1$ by the above discussion,
\item $a_i\le a_{i+1}\le a_i+1$ by Lemma~\ref{thm:head}, and
\item if $a_i=a_{i+1}$, then $a_i=a_{i+1}=\cdots=a_\ell$
    since $\tau^i(\bar B_i) \tau^{i+1}(\bar B_{i+1})$ is left-greedy.
\end{enumerate}
If $a_{\ell-1}\ge\ell$, then
$$(B_{\ell-1}\cdots B_1) P \tau^{-u}\left( B_{\ell-1}\cdots
B_1\right)^{-1}\ge \bD,$$
which contradicts the minimality of $|R|$. So $a_{\ell-1}\le \ell-1$.
Then $a_i=a_{i+1}$ for some $i\le \ell-2$ and so $a_i=a_{i+1}=\cdots =a_\ell\le\ell-1$
so that $\inf(RP\tau^{-u}(R^{-1}))\le 0$. However, by our
choice of $R$, we know that $\inf(RP\tau^{-u}(R^{-1}))> 0$.  Retracing our
steps we conclude that the assumption $H_1 = H_0$ is impossible, so
$H_1 < H_0 < \bD$.

It remains to attack the cases $k>0$.  The method is identical
to the case $k=0$.  Set $V=\cycl^k(W)=\bD^u A_kP_k$ and let $V$
play the role of $W=\bD^u A_0P_0$.
\end{proof}

\smallskip

The proof of Theorem~\ref{thm:cyclbd} is almost complete.  We have learned that
$e< H_m < H_{m-1}< \cdots < H_1 < H_0 < \bD$. This implies that:
 $$0< |H_m| < |H_{m-1}| < \cdots < |H_1| < |H_0| < |\bD|.$$
Thus the length $m+1$ of the chain must be smaller than $|\bD|$, that is
$m+1\le |\bD|-1$. The proof of Theorem~\ref{thm:cyclbd} is complete. \qed

\bigskip

\noindent {\bf Proof of Corollary \ref{cor:complexity estimates}:} The proof
follows directly from Theorem \ref{thm:cyclbd} and the estimates in
\cite{BKL1998}. In Theorem 4.4 of \cite{BKL1998} it is shown that for the new
presentation there is an algorithm rewriting a word into its left greedy form
that is a $O(|W|^2n)$ solution to the word problem. The initial 
preparation of our algorithm puts a given word $W$ into its left
greedy form and takes $O(|W|^2n)$. Notice that the number of
factors is proportional to $|W|$ in the worst case.
In order to compute inf we need to cycle at most $n-2$ times.
After each cycling the new word so-obtained must be put into left greedy
form but this time it takes only $O(|W|n)$ by Corollary 3.14 of
\cite{BKL1998}. Thus the test to determine whether inf is maximal
takes $O(|W|n^2)$. If it is not the entire process must be repeated, but
the number of such repeats, i.e., the total increase of inf, is clearly
bounded by the number of factors so the entire calculation is $O(|W|^2n^2)$.
We note that if $W$ is a positive word, the total increase of inf
is the maximum number of
powers of $\delta$ formed cyclically from $W$ but this number is
clearly bounded by $|W|/(n-1)$, so the entire calculation is
$O(|W|^2n)$.  The discussion for the old presentation is similar and
is left to the reader. \qed

\bigskip

\noindent {\bf Proof of Corollary \ref{cor:geodesic length}:} Let
$W = \bD^u A_1A_2\cdots A_k$ be a word which is in normal form and
which realizes the maximum value $u$ of inf and the minimum value $k$ of sup for
the word class $\{W\}$. The geodesic length $l_Q(\{W\})$ of $\{W\}$ is computed
in
\cite{Charney} (or see \cite{Xu}) as follows:
\begin{itemize}
\item [(i)]If $u\geq 0$ then $W$ is a positive word of geodesic length
$l_Q(\{W\})=u+k$.
\item [(ii)]If $-k\leq u < 0$, then we may use the fact that for every
$X_i\in\bQ$ there exists $Y_i\in\bQ$ with $X_iY_i = \bD$. From this it follows
that $\bD^{-1}X_i = Y_i^{-1}$. Using the additional fact that if $\tau$ is the
index shift automorphism of $\S$2.2, then  $\tau(\bQ) = \bQ$, it follows that we
may eliminate all of the powers of $\bD$ and replace $u$ of the factors
$A_1,A_2,\dots A_u\in\bQ$ with appropriate elements of
$\bQ^{-1}$, thereby achieving a shorter word.  So in this case $l_Q(\{W\}) = k$.
\item [(iii)] If $u<-k$ then every factor $A_1,A_2,\dots A_u\in\bQ$ is replaced
by an appropriate element of $\bQ^{-1}$. After all of these reductions the new
word will be entirely negative. Its geodesic length is $l_Q(\{W\}) = -u$.
\item [(iv)] The three cases may be combined into a single formula: \\
$l_Q(\{W\}) = max(k+u,-u,k).$
\end{itemize}
The above considerations relate to the length of a word class $\{W\}$. However,
observe  that the normal  form for elements in the conjugacy
class $[W]$ is identical to that for the word class, moreover if $Y,Z$ are in the
super summit set of $[W]$ then $\inf(Y) = \inf(Z)$ and $\sup(Y) = \sup(Z)$.  Since
the complexity of computing
$l_Q([W])$ is identical to the complexity of computing $\inf([W])$ and
$\sup([W])$, the assertion then follows from Corollary \ref{cor:complexity
estimates}. \qed

\section{{\bf Are the cycling-decycling bounds sharp?}}
\label{sec:ex}
Note that the bound we obtained for the number of cyclings and decyclings
in Theorem~\ref{thm:cyclbd} is $n-2$ for the new presentation and 
$-1 + (n-1)(n-2)/2$ for the old presentation. In this section we investigate
whether these bounds are sharp.

\begin{figure}
$$  {\includegraphics{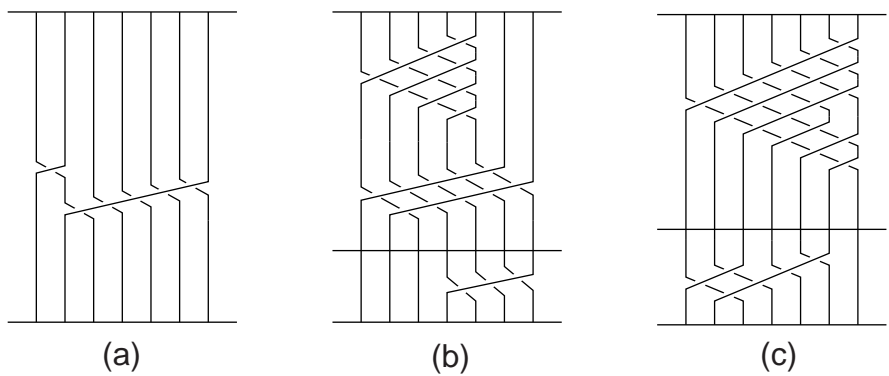}} $$
\caption{Examples}\label{fig:ex}
\end{figure}

 We first give an example of $n$-braid written in the
new generators for which $n-2$ cyclings are required to increase the infimum.
This shows that the bound given in Theorem~\ref{thm:cyclbd} is sharp for
the new presentation. To simplify notation, use $[t, t-1,\ldots,
s]$ instead of $a_{t(t-1)} a_{(t-1)(t-2)}\cdots a_{(s+1)s}.$
Consider the example $W= ([2,1][5,4,3])([3,2])$ in normal form. Then

\begin{eqnarray*}
\cycl(W) &=& ([3, 2]) ([2,1][5,4,3]) = ([3,2,1][5,4]) ([4,3])\\
\cycl^2(W) &=& ([4,3]) ([3,2,1][5,4]) = ([4,3,2,1]) ([5,4])\\
\cycl^3(W) &=& ([5,4]) ([4,3,2,1]) = [5,4,3,2,1] = \delta
\end{eqnarray*}
So $\inf(W)=\inf(\cycl(W)) = \inf(\cycl^2(W))=0$ but $\inf(\cycl^3(W))=1$.
More generally, if
$$W=[2,1] [n,n-1,\ldots,3,2] = ([2,1] [n,n-1,\ldots, 3])([3,2]),$$
then $\inf(W)=\inf(\cycl^{n-3}(W))=0$ but $\inf(\cycl^{n-2}(W))=1$.
See Figure~\ref{fig:ex}(a) for a sketch of the braid $W$ in the case $n=7$.

In the old presentation, the example in~\cite{EM1994}
shows that $\inf(W)=\inf(\cycl(W))=0$ but $\inf(\cycl^2(W))=1$.
There are plenty of examples for which more than 2 cyclings are required
to increase the infimum.
Let $(a_1,\ldots, a_n)$ denote the permutation braid corresponding to
the permutation $\pi$ on $\{1,\ldots,n\}$ defined by $\pi(i)=a_i$,
Consider the following example, with $W\in B_{2k+1}$.
$$W = (\underbrace{2k+1, 2k, \ldots , 3}_{2k-1}, 1, 2)
   (\underbrace{1, 2, \ldots , k}_k, \underbrace{k+2, \ldots , 2k+1}_k, k+1)$$
Then $\inf(W)=\inf(\cycl^{2k-1}(W))=0$ but $\inf(\cycl^{2k}(W))=1$.
See Figure \ref{fig:ex}(b) for a sketch of this example in the case $n=7$.

For another example let $W\in B_{2k+1}$ be such that
\begin{eqnarray*}
W&=&(\underbrace{2k+1,\ldots,k+3}_{k-1}, k+1, k+2,\underbrace{k, k-1, \ldots, 1}_k)\\
 &&\qquad (\underbrace{3,4,\ldots,k+1}_{k-1}, 1, \underbrace{k+2,\ldots,2k}_{k-1},2,2k+1)
\end{eqnarray*}
See Figure~\ref{fig:ex}(c). Then $4k-5$ cyclings are needed to
increase the infimum. So if $n$ is odd, there is an example for
which $2n-7$ cyclings are needed. Therefore the lower bound for
the old presentation is at least linear in $n$.

We do not know an exact bound that works for every $n$-braid
written in the old generators. It is easy to see that the upper
bound in $B_3$ is 1, that is, if $\inf(W)=\inf(\cycl(W))$ for $W\in B_3$,
then the infimum is already maximized. In an exhaustive search, we
learned that the upper bound for $B_4$, using the old presentation,
is 2 for positive words whose normal form contains up to 5 canonical factors.

\section{{\bf Complexity issues and the conjugacy
problem}}
\label{sec:complexity issues}

In this section we
consider implications of the work in the preceding sections for the complexity
of the conjugacy problem in $B_n$.

\subsection{The special cases $n=3$ and $4$:}

 Before discussing the
problem, it will be helpful to review what is known about the cases $n=3$ and
$4$, since well-chosen examples always help one to arrive at a better
understanding of a problem.  In the manuscript
\cite{Xu} P.J. Xu introduced the new presentation for $B_3$ and used it
to solve the word and conjugacy problems in $B_3$ and to study the
letter lengths of shortest words in a word and conjugacy class in $B_3$, using
the new presentation.  Her main result in this regard was that words of shortest
`geodesic length' (she doesn't use the term geodesic length, which was
introduced after she completed her work)  are, without further work, also words
of shortest letter length in the new generators. She also found growth functions
for $B_3$, both for word classes and conjugacy classes, proved that they were
rational,  and computed the rational functions which described them.  Her
algorithm for the conjugacy problem was clearly polynomial in $|W|$.

In \cite{KKL} the word and conjugacy problems were solved in $B_4$, using the
new presentation and following the methods of \cite{Xu}. The authors also
solved the shortest word problem in conjugacy classes. In a forthcoming paper
the second and third author of this paper will prove that the algorithm for the
conjugacy problem in \cite{KKL} is polynomial in $|W|$.

\smallskip

\subsection{The conjugacy problem:}  In $\S$2.10 above the super
summit set SSS$([W])$ of the conjugacy class of $W\in B_n$ is defined. It is a
finite set and it can be computed in a systematic manner in a finite number of
steps from any braid word
$W$ which realizes $\inf([W])$ and $\sup([W])$. The Theorem which is quoted in
$\S$2.11 above asserts that
$W$ is conjugate to $V$ in $B_n$ if and only if
$\inf([W]) = \inf([V]), \ \ \sup([W]) = \sup([V])$ and SSS$([W]) =
{\rm SSS}([V])$.

The super summit set has a fairly transparent structure when $n=3$, the main
reason being that words in $\bQ^\star = \bQ\setminus \{\delta,e\}$ all have
length 1. In $B_4$ the situation is a little bit more complicated, but still
within reach.
Let $W = \delta^u A_1A_2\cdots A_k$ be in the super summit set of $[W]$ and be in
normal form, and let $A = A_1A_2\cdots A_k$ be the `positive part' of $W$.
Notice that $\delta$ has letter length 3 and the $A_j's$ are elements of
$\bQ^\star$, and so have letter length 1 or 2. Let $k_1$ (resp. $k_2$) be the
number of factors in $A$ which have length 1 (resp. 2). Let
$e$ be the exponent sum of $W$. Clearly $e$ is a class invariant.  Since $k =
k_1 + k_2$ and since $e = 3u + k_1 + 2k_2$ it follows that $k_1$ and $k_2$ are
determined by the triplet $(u=\inf, k=\sup, e)$.  This makes the SSS
somewhat easier to understand in the case $n=4$ than in the general case.

In the general case the super summit set SSS$([W])$ splits into orbits under
cycling and decycling. Clearly the number of such orbits and their sizes are
class invariants, but unfortunately we have examples to show that they are not
complete invariants. The orbits are complicated by the fact that $\bQ^\star$
contains elements of letter length $1,2,\dots,n-2$, and the number $k_i$ of
elements of letter length $i$  of a
member of SSS$([W])$ is no longer controlled by $(u,k,e)$. We don't know
whether $k_1,\dots,k_{n-2}$ are orbit invariants, and if they can vary from one
orbit to another. Also, while it is known that one can pass from any orbit to
any other orbit by conjugating by an appropriate product of elements of
$\bQ$, it is difficult to understand which products do the job. While the
super summit set is a great improvement over the summit set of
\cite{Garside}, it is still too big to make it possible to study many
examples.  For all these reasons the complexity of the conjugacy problem remains
open at this time.  Nevertheless, based on what we know, we conjecture:
\begin{con}
\label{con:complexity of the conjugacy problem}
There is an algorithmic solution to the conjugacy problem in $B_n$, using the
combinatorial approach which is described in this paper,  which is polynomial in
word length $|W|$ for each fixed braid index $n$.
\end{con}

\end{document}